\documentclass[12pt,a4paper]{article}
\usepackage{amsmath,amstext,amssymb,amscd}

\usepackage{amsmath}
\usepackage{mathtools}
\usepackage{breqn}

\newcommand{\Xcomment}[1]{}

\begin{document}

\baselineskip=15pt
\parskip=2pt

 \title{Erratum to ``$B_2$-crystals: Axioms, structure, models'' \\
 \lbrack J. Combin. Theory, Ser. A, 116 (2009), 265--289\rbrack}

 \author{Vladimir I.~Danilov
\thanks{Central Institute of Economics and
Mathematics of the RAS, 47, Nakhimovskii Prospect, 117418 Moscow, Russia;
email: danilov@cemi.rssi.ru.}
 \and
Alexander V.~Karzanov
\thanks{Central Institute of Economics and Mathematics of
the RAS, 47, Nakhimovskii Prospect, 117418 Moscow, Russia; email:
akarzanov7@gmail.com. Corresponding author. }
  \and
Gleb A.~Koshevoy
\thanks{Institute for Information Transmission Problems of
the RAS, 19, Bol'shoi Karetnyi per., 127051 Moscow, Russia; email:
koshevoyga@gmail.com. }
  }

\date{}

 \maketitle

 \begin{quote}
 {\bf Abstract.} In this erratum we explain how to implement two axioms  stated in~\cite{DKK} so as to get a purely ``local'' characterization for finite $B_2$-crystals, which was declared but not clarified at some moments there. Also we correct some inaccuracies in that paper.
 \end{quote}


\section{Local axioms}  \label{sec:axioms}

Crystal graphs of type $B_2$ constitute the simplest case of doubly laced Kashiwara's crystals in representation theory, and paper~\cite{DKK} was devoted to a combinatorial study of this class of 2-edge-colored directed graphs. Before the appearance of that paper, there have been known several ``global'' characterizations of $B_2$-crystals (e.g., using generalized Young tableaux, Lusztig's canonical bases,  Littelmann's path model), and one of the main purposes of~\cite{DKK} was to find a characterization in ``local terms''. Here we understand by the ``locality'' a way of defining a crystal graph only via requirements on the structure of small neighborhoods of vertices (in the sense that the radius and size of a neighborhood are bounded by a constant).

Theorem~4 in~\cite{DKK} characterizes the set of so-called S-graphs by axioms (B0)--(B4), (B$'$3),  (B$'$4) and (BA) (where the last axiom reduces to a series of other ones), which are then shown in the paper to be precisely the set of $B_2$-crystals, in both finite and infinite cases. It is seen that a majority of these axioms are obviously local, but there are two axioms, namely, (B1) and (B2), for which a possibility to be implemented by using merely local terms is not clarified in~\cite{DKK}. 

Next we explain how to fulfill this task in the assumption that the input graphs are finite and acyclic (i.e., having no directed cycles).

More precisely, we deal with a finite acyclic directed graph $G=(V,E)$ in which the edge set $E$ is partitioned into two subsets $E_1$ and $E_2$, consisting of edges of \emph{color 1} and \emph{color 2}, also called  \emph{1-edges} and \emph{2-edges}, respectively. Also it is usually assumed that $G$ is weakly connected. According to Axiom~(B0), for $i=1,2$, each vertex $v\in V$ has at most one entering and at most one leaving edge of color $i$. Therefore, the subgraph $(V,E_i)$ consists of pairwise disjoint directed paths covering all vertices, called $i$-\emph{strings}. Besides, some vertices and edges are distinguished as \emph{central} ones. Then Axioms (B1) and (B2) read as follows.

 \begin{itemize}
  \item[(B1)] Each 1-string has exactly one central element, which is either a vertex or an edge.
  \end{itemize}

This partitions $V$ into three subsets, consisting of central, left and right vertices, where a vertex is called \emph{left} (resp. \emph{right}) if it lies in its 1-string before (resp. after) the central element. (Note that when $(u,v)$ is a central edge, its beginning vertex $u$ is regarded as left, while the end vertex $v$ as right vertex.) Accordingly, a non-central 1-edge $(u,v)$ is called \emph{left} (\emph{right}) if $u$ is left (resp. $v$ is right). 

\begin{itemize}
\item[(B2)]
Each 2-string $P$ contains exactly one central vertex $v$. Moreover,
all vertices of $P$ lying {\em before} $v$ are right, whereas all
vertices lying {\em after} $v$ are left.
  \end{itemize}
  
Formally speaking, Axioms (B1) and (B2) are not local. In order to obtain their local implementations, we assume that an input graph $G=(V,E)$ as before is equipped with \emph{labels} $\ell(v)$ on the vertices $v\in V$ which take values in the 3-element set $\{0,c,1\}$. We impose the following local requirements on $(G,\ell)$.

 \begin{itemize}
  \item[(B1(i))]  For each 1-edge $(u,v)$, the pair $(\ell(u),\ell(v))$ is equal to one of $(0,0),(0,c),(0,1)$, $(c,1),(1,1)$.
  
  \item[(B1(ii))] If a vertex $v$ has no entering 1-edge, then $\ell(v)\ne 1$ (i.e., $\ell(v)\in\{0,c\}$); and if a vertex $v$ has no leaving 1-edge, then $\ell(v)\ne 0$.
  
  \item[(B2(i))] For each 2-edge $(u,v)$, the pair $(\ell(u),\ell(v))$ is equal to one of $(1,1),(1,c)$, $(c,0),(0,0)$.  
  
  \item[(B2(ii))] If a vertex $v$ has no entering 2-edge, then $\ell(v)\ne 0$; and if a vertex $v$ has no leaving 2-edge, then $\ell(v)\ne 1$.  
  \end{itemize}

\noindent In particular, (B1(ii)) implies that if $v$ has neither entering nor leaving 1-edge, then $\ell(v)=c$, and similarly for (B2(ii)). (Note that assigning labels as above looks in spirit of the local axiomatics for simply laced crystals given in~\cite{stem} where 0,1 labels are assigned to the edges of a crystal graph.)

The vertices labeled $0,c,1$ are naturally interpreted as left, central and right ones, respectively. Accordingly, a 1-edge $(u,v)$ is regarded as left if $(\ell(u),\ell(v))\in\{(0,0),(0,c)\}$, central if $(\ell(u),\ell(v))=(0,1)$, and right if $(\ell(u),\ell(v))\in\{(c,1),(1,1)\}$. As to a 2-edge $(u,v)$, it is regarded as left if $(\ell(u),\ell(v))\in\{(c,0),(0,0)\}$, and right if $(\ell(u),\ell(v))\in\{(1,1),(1,c)\}$.  

As an easy consequence of the above assignments, we conclude with the following
 \medskip
 
\noindent \textbf{Proposition} ~\emph{For finite acyclic graphs satisfying~(B0), Axioms (B1)--(B2) are equivalent to imposing (B1(i),(ii)),(B2(i),(ii)).}
 \medskip

Thus, when dealing with finite acyclic graphs, we obtain a ``purely local'' axiomatics for finite $B_2$-crystals, as required. Note that the formal requirement that an input graph $G$ is acyclic can be realized by imposing additional local variables and constraints. Namely, let us endow each vertex $v$ with an additional number $\pi(v)$ (a ``potential'') and impose the condition: $\pi(u)<\pi(v)$ for each edge $(u,v)$ of $G$. Clearly $G$ is acyclic if and only if a feasible $\pi$ does exists.

However, the above labeling method does not work when an input graph is infinite (since in this case Axioms (B1(i),(ii)),(B2(i),(ii)) do not forbid the existence of infinite monochromatic strings without central elements). Therefore, our local characterization is applicable only to finite $B_2$-crystals.

We finish this section with one more useful observation.
 \medskip
 
 \noindent\textbf{Remark.}  In fact,~\cite{DKK} exhibits two local axiomatics for $B_2$-crystals. One of them is just what was discussed above. An alternative local axiomatics is provided by the so-called \emph{worm model} developed in Section~4 of~\cite{DKK}. In this model, each vertex $v$ of a graph is endowed with a six-tuple $\tau(v)$ of integers $(x,x',x'',y,y',y'')$ satisfying conditions ~(A),(B),(C) on page~278 (due to this, one can associate $v$ with a figure consisting of at most two line segments (horizontal and vertical ones) in a rectangle of the plane). Local conditions in (i)--(vi) (on pages 278--279) prescribe how the six-tuple $\tau(v)$ can change under the action of crystal operators (thus defining the edges of colors 1 and 2 entering and leaving $v$).


\section{Other corrections}  \label{sec:other}

In this section we correct three inaccuracies from~\cite{DKK}; they were pointed out by Shunsuke Tsuchioka in~\cite{shu}.

The first one concerns Remark~2 on page 273. It gives a commentary to axiom (B3) and Corollary~2 and has incorrect sentences in the second paragraph, which can be disproved by elementary examples. Nevertheless, this remark does not affect the main content of the paper; moreover, there is no statement in the remaining part of the paper where elements from this remark are quoted or used. Due to this, the second paragraph in Remark~2 should be deleted, which causes no flaw in the whole content. 

The second and third inaccuracies involve Corollaries~2 and~3 on page~273 in~\cite{DKK}. They read as follows.
\medskip

\noindent \textbf{Corollary 2}  \emph{
~Let $(u,v)$ be a central 1-edge. Then there are a 2-edge $(u',u)$ and a
2-edge $(v,v')$. Moreover, both vertices $u',v'$ are central. }
 \medskip

\noindent \textbf{Corollary 3}  
\emph{Let $(u,v)$ be a central 1-edge. Let a 2-edge $(u,w)$ leave $u$. Then
there is a 1-edge $(w,w')$, and the vertex $w'$ is central.}
\medskip

Both corollaries are correct but their proofs given in~\cite{DKK} contain gaps. Correct proofs essentially use Axiom (B4) (and its dual (B$'$4)), and accordingly, these corollaries should be placed after this axiom. Below we give correct proofs.
  \medskip
  
\noindent\textbf{Proof of Corollary~2}
~Since the edge $(u,v)$ is central, the vertex $u$ is left and the vertex $v$ is right. Hence $u$ has an entering 2-edge, $(u',u)$ say, and $v$ has a leaving 2-edge, $(v,v')$ say; note that the former edge is left and the latter is right. Suppose that the vertex $u'$ is not central. Then $u'$ is left, and therefore $u'$ has a leaving 1-edge, $(u',u'')$ say. The edge $(u',u'')$ cannot be left; for otherwise $u',u'',u,v$ would give a commutative square (by Axiom (B3)), whence $(u,v)$ should be left, not central. Therefore, $(u',u'')$ must be central. But then, applying Axiom (B4) to the edges $(u',u''),(u',u),(u,v),(v,v')$, we conclude that the vertex $v$ cannot be right, contradicting the condition that $(u,v)$ is central. Thus, $u'$ is central, as required.
 
The assertion that $v'$ is central is symmetric.
 \medskip
  
\noindent\textbf{Proof of Corollary~3} ~Since $u$ is left, $w$ is left as well. So a 1-edge $(w,w')$ does exist. This edge cannot be central; for otherwise the vertex $u$ would be central, by Corollary~2, contradicting the condition that $(u,v)$ is central. Suppose that $w'$ is not central. Then $w'$ is left, and therefore, $w'$ has an entering 2-edge, $(w'',w')$ say. This edge is left and we can apply~(B3)(ii) to the edges $(w'',w')$ and $(w,w')$, obtaining a commutative square containing the vertices $w,w',w''$. This square must contain the vertices $u,v$ as well, implying $v=w''$ and leading to a contradiction to the fact that the edge $(u,v)$ is central.


\end{document}